\RequirePackage{fix-cm}
\documentclass{svjour3}                     

\smartqed  

\usepackage{graphicx, booktabs}
\usepackage{cite}

\usepackage{mathtools,amsfonts,dsfont}

\let\ALPHABET\mathcal
\let\VEC\mathbf

\def\DEFINED{\coloneqq}
\def\EXP{\mathds{E}}

\def\PR {\mathds{P}}

\newcommand\naturalnumbers{\mathbb{Z}^{+}}
\newcommand\reals{\mathbb{R}}

\begin{document}

\title{Decentralized stochastic control%
  \thanks{This work was supported by Fonds de recherche du Qu\'ebec--Nature et
  technologies (FRQ-NT) Establishment of New Researcher Grant 166065 and by
  Natural Science and Engineering Research Council of Canada (NSERC) Discovery
  Grant 402753-11.}}


\author{Aditya Mahajan \and Mehnaz Mannan}


\institute{A. Mahajan and M. Mannan \at
           3480 University, ECE Department, McGill University, Montreal, QC, Canada, H3A 0E9 \\
           E-mail: aditya.mahajan@mcgill.ca, mehnaz.mannan@mail.mcgill.ca
}


\maketitle

\begin{abstract}
  Decentralized stochastic control refers to the multi-stage optimization of a
  dynamical system by multiple controllers that have access to different
  information. Decentralization of information gives rise to new conceptual
  challenges that require new solution approaches. In this expository paper, we
  use the notion of an \emph{information-state} to explain the two commonly used
  solution approaches to decentralized control: the person-by-person approach
  and the common-information approach.

\end{abstract}

\section {Introduction} \label{sec:intro}

Centralized stochastic control refers to the multi-stage optimization of a
dynamical system by a single controller. Stochastic control, and the associated
principle of dynamic programming, have roots in statistical sequential
analysis~\cite{ArrowBlackwellGirshick:1949} and have been used in various
application domains including operations research~\cite{Powell:2007},
economics~\cite{Stokey:1989}, engineering~\cite{Bertsekas:1995}, computer
science~\cite{RussellNorvig:1995}, and mathematics~\cite{Bellman:1957}.  The
fundamental assumption of centralized stochastic control is that the decisions
at each stage are made by a single controller that has \emph{perfect recall},
that is, a controller that remembers its past observations and decisions. This
fundamental assumption is violated in many modern applications where decisions
are made by multiple controllers. The multi-stage optimization of such
systems is called \emph{decentralized  stochastic control}.

Decentralized stochastic control started with seminal work of Marschak and
Radner~\cite{MarschakRadner:1972, Radner:1962} on static systems that arise in
organizations and of Witsenhausen~\cite{Witsenhausen:1971, Witsenhausen:1971a,
Witsenhausen:1973} on dynamic systems that arise in systems and control. We
refer the reader to~\cite{BasarBansal:1989,Ho:1980} for a discussion of the
history of decentralized stochastic control and
to~\cite{YukselBasar:2013,MMRY:tutorial-CDC,Oliehoek:JAIR2013} for survey of
recent results. 

Decentralized stochastic control is fundamentally different from and
significantly more challenging than centralized stochastic control. Dynamic
programming, which is the primary solution concept of centralized
stochastic control, does not directly work in decentralized stochastic control
and new ways of thinking need to be developed to address information
decentralization. The focus of this expository paper is to highlight the
conceptual challenges of decentralized control and explain the intuition behind the solution
approaches. No new results are presented in this paper; rather we present new
insights and connections between existing results. Since the focus is on
conceptual understanding, we do not present proofs and ignore the
technical details, in particular, measurability concerns, in our description.

We use the following notation. Random variables are denoted by upper case
letters; their realizations by the corresponding lower case letters; and their
space of realizations by the corresponding calligraphic letters. For integers $a
\le b$, $X_{a:b}$ is a short hand for the set $\{X_a, X_{a+1}, \dots, X_b\}$.
When $a > b$, $X_{a:b}$ refers to the empty set. In general, subscripts are used
as time index while superscripts are used to index controllers. $\PR(\cdot)$
denotes the probability of an event and $\EXP[\cdot]$ denotes the expectation of
a random variable. For a collection of functions $\VEC g$, the notations
$\PR^{\VEC g}(\cdot)$ and $\EXP^{\VEC g}[\cdot]$ indicate that the probability
measure and the expectation depend on the choice of the functions $\VEC g$.
$\naturalnumbers$ denotes the set of positive integers and
$\reals$ denotes the set of real numbers. 

\section {Decentralized stochastic control: Models and problem formulation}
\label{sec:decentralized}

\subsection{State, observation, and control processes}
Consider a dynamical system with $n$~controllers.
Let $\{X_t\}_{t=0}^\infty$, $X_t \in \ALPHABET X$, denote the state process of
the system. Controller ~$i$, $i \in
\{1,\dots,n\}$, causally observes the process $\{Y^i_t\}_{t=0}^\infty$, $Y^i_t
\in \ALPHABET Y^i$, and generates a control process $\{U^i_t\}_{t=0}^\infty$,
$U^i_t \in \ALPHABET U^i$. The system yields a reward $\{R_t\}_{t=0}^\infty$.
These processes are related as follows:
\begin{enumerate}
  \item Let $\VEC U_t \DEFINED \{U^1_t, \dots, U^n_t\}$ denote the control action of
    all controllers at time~$t$. Then, the reward at time~$t$ depends only on
    the current state $X_t$, the future state $X_{t+1}$, and the current control
    actions $\VEC U_t$. Furthermore, the state process $\{X_t\}_{t=0}^\infty$ is
    a controlled Markov process given $\{\VEC U_t\}_{t=0}^\infty$, i.e., for any
    $\ALPHABET A \subseteq \ALPHABET X$ and $\ALPHABET B \subseteq \mathbb R$,
    and any realization $x_{1:t}$ of $X_{1:t}$ and $\VEC u_{1:t}$ of $\VEC
    U_{1:t}$, we have that
    \begin{multline}
      \label{eq:d-dynamics}
      \PR(X_{t+1} \in \ALPHABET A, R_t \in \ALPHABET B \mid X_{1:t} = x_{1:t}, \VEC
      U_{1:t} = \VEC u_{1:t})
      \\ =
      \PR(X_{t+1} \in \ALPHABET A, R_t \in \ALPHABET B \mid X_{t} = x_{t}, \VEC
      U_{t} = \VEC u_{t}).
    \end{multline}

  \item The observations $\VEC Y_t \DEFINED \{Y^1_t, \dots, Y^n_t\}$ depend only on
    current state $X_t$ and previous control actions $\VEC U_{t-1}$, i.e., 
    for any $\ALPHABET A^i \subseteq \ALPHABET Y^i$ and any realization $x_{1:t}$ of
    $X_{1:t}$ and $\VEC u_{1:t-1}$ of $\VEC U_{1:t-1}$, we have that
    \begin{multline}
      \label{eq:d-observation}
      \PR\Big(\VEC Y_t \in \prod_{i=1}^n \ALPHABET A^i \Bigm| X_{1:t} = x_{1:t}, \VEC
      U_{1:t-1} = \VEC u_{1:t-1}\Big)
      \\ =
      \PR\Big(\VEC Y_{t} \in \prod_{i=1}^n \ALPHABET A^i \Bigm|
          X_{t} = x_{t}, \VEC U_{t-1} = \VEC u_{t-1}\Big).
    \end{multline}
\end{enumerate}

\subsection {Information structure}
\label{sec:info-structure}

At time~$t$, controller~$i$, $i \in \{1,\dots,n\}$, has access to information
$I^i_t$ which is a superset of the history $\{Y^i_{1:t}, U^i_{1:t-1}\}$ of the
observations and control actions at controller~$i$ and a subset of the history
$\{\VEC Y_{1:t}, \VEC U_{1:t-1}\}$ of the observations and control actions at
all controllers, i.e., 
\[
  \{Y^i_{1:t}, U^i_{1:t-1}\} \subseteq I^i_t \subseteq
  \{\VEC Y_{1:t}, \VEC U_{1:t-1}\}.
\]
The collection $(I^i_t, i \in \{1,\dots,n\}, t = 0,1,\dots)$, which  is called the
\emph{information structure} of the system, captures
\emph{who knows what about the system and when}. A decentralized system is characterized by its
information structure. Some examples of information structures are given below.
For ease of exposition, we use $J^i_t$ to denote $\{Y^i_{1:t}, U^i_{1:t-1}\}$
and refer to it as \emph{self information}. 
\begin{enumerate}
  \item \emph{Complete information sharing information structure} refers to a
    system in which each controller has access to the self information of all
    other controllers, i.e.,
    \[
      I^i_t = \bigcup_{j=1}^n J^j_t, \quad \forall i \in \{1,\dots,n\}.
    \]

  \item \emph{$k$-step delayed sharing information structure} refers to a system
    in which each controller has access to $k$-step delayed self information of
    all other controllers, i.e.,
    \[
      I^i_t = J^i_t \cup \Big( 
      \bigcup_{\substack{ j = 1 \\ j \neq i }}^n J^j_{t-k}
      \Big), \quad \forall i \in \{1, \dots, n\}.
    \]

  \item \emph{$k$-step periodic sharing information structure} refers to a system
    in which all controllers periodically share their self information after
    every $k$ steps, i.e.,
    \[
      I^i_t = J^i_t \cup \Big( 
      \bigcup_{\substack{ j = 1 \\ j \neq i }}^n J^j_{\lfloor t/k \rfloor k}
      \Big), \quad \forall i \in \{1, \dots, n\}.
    \]

  \item \emph{No sharing information structure} refers to a system in which the
    controllers do not share their self information, i.e.,
    \[
      I^i_t = J^i_t, \quad \forall i \in \{1, \dots, n\}.
    \]
\end{enumerate}

\subsection{Control strategies and problem formulation}

Based on the information $I^i_t$ available to it, controller~$i$ chooses action
$U^i_t$ using a \emph{control law} $g^i_t \colon I^i_t \mapsto U^i_t$. The
collection of control laws $\VEC g^i \DEFINED (g^i_0, g^i_1, \dots)$ is called a
\emph{control strategy of controller~$i$}. The collection $\VEC g \DEFINED (\VEC g^1,
\dots, \VEC g^n)$ is called the \emph{control strategy of the system}.

The optimization objective is to pick a control strategy $\VEC g$ to maximize the expected
discounted reward
\begin{equation}
  \label{eq:d-reward}
  \Lambda(\VEC g) \DEFINED \EXP^{\VEC g}\Big[ \sum_{t=0}^\infty \beta^t R_t \Big]
\end{equation}
for a given discount factor $\beta \in (0,1)$.

\subsection {An example}
\label{sec:d-example}

To illustrate these concepts, let's consider a stylized example of a
communication system in which two devices transmit over a multiple access
channel.

\begin{description}
  \item[\emph{Packet arrival at the devices.}]
    Packets arrive at device~$i$, $i \in \{1,2\}$, according to Bernoulli
    processes $\{W^i_t\}_{t=0}^\infty$ with success probability $p^i$.
    Device~$i$ may store $N^i_t \in \{0,1\}$ packets in a buffer. If a packet
    arrives when the buffer is full, the packet is dropped.

  \item[\emph{Channel model.}]
    At time $t$, the channel-state $S_t \in \{0,1\}$ may be idle ($S_t = 0$) or
    busy ($S_t = 1$). The channel-state process $\{S_t\}_{t=0}^\infty$ is a
    Markov process with known initial distribution and transition matrix
    $\mathbf P = \left[\begin{smallmatrix}
      \alpha_0 & 1 - \alpha_0 \\
      1 - \alpha_1 & \alpha_1
    \end{smallmatrix}\right]$. The
    channel-state process is independent of the packet-arrival process at the
    device.

  \item[\emph{System dynamics.}]
    At time~$t$, device~$i$, $i \in \{1,2\}$, may transmit $U^i_t \in \{0,1\}$
    packets, $U^i_t \le N^i_t$. If only one device transmits and the channel is
    idle, the transmission is successful and the transmitted packet is removed
    from the buffer. Otherwise the transmission is unsuccessful. The state of
    each buffer evolves as
    \begin{equation} \label{eq:d-buffer}
      N^i_{t+1} = \min\{N^i_t - U^i_t(1 - U^j_t) (1 - S_t) + W^i_t, 1\},
      \quad \forall i \in \{1,2\}, \quad j = 3 - i.
    \end{equation}
    Each transmission costs $c$ and a successful transmission yields a reward
    $r$. Thus, the total reward \emph{for both devices} is
    \[
      R_t = - (U^1_t + U^2_t)c + (U^1_t \oplus U^2_t)(1-S_t) r
    \]
    where $\oplus$ denotes the XOR operation.

  \item[\emph{Observation model.}]
    Controller~$i$, $i \in \{1,2\}$, perfectly observes the number $N^i_t$ of
    packets in the buffer. In addition, \emph{both} controllers observe the
    one-step delayed control actions $(U^1_{t-1}, U^2_{t-1})$ of each other and
    the channel state if \emph{either of devices transmit}. Let $H_t$ denote
    this additional observation. Then $H_t = S_{t-1}$ if $U^1_{t-1} + U^2_{t-1}
    > 0$, otherwise $H_t = \mathfrak E$ (which denotes no channel-state
    observation).

  \item[\emph{Information structure and objective.}]
    The information $I^i_t$ available at device $i$, $i \in \{0,1\}$, is given
    by $I^i_t = \{N^i_{1:t}, H_{1:t},\allowbreak U^1_{1:t-1}, U^2_{1:t-1}\}$.
    Based on the information available to it, device~$i$ chooses control action
    $U^i_t$ using a control law $g^i_t \colon I^i_t \mapsto U^i_t$. The
    collection of control laws $(\VEC g^1, \VEC g^2)$, where $\VEC g^i \DEFINED (g^i_0,
    g^i_1, \dots)$, is called a \emph{control strategy}. The objective is to
    pick a control strategy $(\VEC g^1, \VEC g^2)$ to maximize the expected
    discounted reward
    \[
      \Lambda(\VEC g^1, \VEC g^2) \DEFINED
      \EXP^{(\VEC g^1, \VEC g^2)}\Big[ \sum_{t=0}^\infty \beta^t R_t \Big].
    \]
\end{description}
We make the following assumption in the paper.
\begin{enumerate}
  \item [\textbf{(A)}] The arrival process at the two controllers is
    independent.
\end{enumerate}

\subsection {Conceptual difficulties in finding an optimal solution}
\label{sec:difficulties}

There are two conceptual difficulties in the optimal design of decentralized
stochastic control:
\begin{enumerate}
  \item The optimal control problem is a functional optimization problem where
    we have to choose an infinite sequence of control laws $\VEC g$ to maximize
    the expected total reward.

  \item In general, the domain $I^i_t$ of control laws $g^i_t$ increases with
    time. Therefore, it is not immediately clear if we can solve the above
    optimization problem; even if it is solved, it is not immediately clear if
    we can implement the optimal solution.
\end{enumerate}

Similar conceptual difficulties arise in centralized stochastic control where
they are resolved by identifying an appropriate \emph{information-state}
process and solving a corresponding dynamic program. It is not possible to
directly apply such an approach to decentralized stochastic control problems. 

In order to better understand the difficulties in extending the solution
techniques of centralized stochastic control to decentralized stochastic
control, we revisit the main results of centralized stochastic control in the
next section. 

%
%
%
%

\section {Overview of centralized stochastic control}
\label{sec:centralized}

A centralized stochastic control system is a special case of a decentralized
stochastic control system in which there is only one controller ($n=1$), and
the controller has perfect recall ($I^1_t \subseteq I^1_{t+1}$), i.e., the
controller remembers everything that it has seen and done in the past. For ease
of notation, we drop the superscript~$1$ and denote the observation,
information, control action, and control law of the controller by $Y_t$, $I_t$,
$U_t$, and $g_t$, respectively. Using this notation, the information available
to the controller at time~$t$ is given by $I_t = \{Y_{1:t}, U_{1:t-1}\}$. The
controller uses a control law $g_t \colon I_t \mapsto U_t$ to choose a control
action $U_t$. The collection $\VEC g = (g_0, g_1, \dots)$ of control laws 
is called a \emph{control strategy}.

The optimization objective is to pick a control strategy $\VEC g$ to maximize the expected
discounted reward
\begin{equation}
  \label{eq:reward}
  \Lambda(\VEC g) \DEFINED \EXP^{\VEC g}\Big[ \sum_{t=0}^\infty \beta^t R_t \Big]
\end{equation}
for a given discount factor $\beta \in (0,1)$.

In the centralized stochastic control literature, the above model is sometimes
referred to a partially observable Markov decision process (POMDP). The
solution to a POMDP is obtained in two steps~\cite{Bertsekas:1995}.
\begin{enumerate}
  \item Consider a simpler model in which the controller perfectly observes the
    state of the system, i.e., $Y_t = X_t$. Such a model is called a Markov
    decision process (MDP). Show that there is no loss of optimality in
    restricting attention to \emph{Markov strategies}, i.e., control laws of the
    form $g_t \colon X_t \mapsto U_t$. Obtain an optimal control strategy of
    this form by solving an appropriate dynamic program.

  \item Define a \emph{belief state} of a POMDP as the posterior distribution of
    $X_t$ given the information at the controller, i.e., $B_t(\cdot) = \PR(X_t =
    \cdot \mid I_t)$. Show that the belief state is a MDP, and use the results
    for MDP.
\end{enumerate}

An alternate (and, in our opinion, a more transparent) approach is identify an
\emph{information-state} process of the system and present the solution in terms
of the information state. We present this approach below.

\begin{definition}
  \label{def:info}
  A process $\{Z_t \}_{t=0}^\infty$, $Z_t \in \ALPHABET Z_t$, is
  called an \emph{information-state} process if it satisfies the following
  properties:
  \begin{enumerate}
    \item $Z_t$ is a function of the information $I_t$ available at time~$t$,
      i.e., there exist a series of functions $\{f_t\}_{t=0}^\infty$ such that
      \begin{equation} \label{eq:info-def}
        Z_t = f_t(I_t).
      \end{equation}

    \item The process $Z_t$ is a controlled Markov process controlled by
      $\{U_t\}_{i=0}^\infty$, that is for any $\ALPHABET A \subseteq \ALPHABET
      Z_{t+1}$ and any realization $i_t$ of $I_t$ and any choice $u_t$ of $U_t$,
      we have that
      \begin{equation}
        \label{eq:info-evolution}
        \PR(Z_{t+1} \in \ALPHABET A \mid I_t = i_t, U_t = u_t) 
        =
        \PR(Z_{t+1} \in \ALPHABET A \mid Z_t = f_t(i_t), U_t = u_t).
      \end{equation}

    \item $Z_t$ absorbs the effect all the available information on the current
      rewards, i.e., for any $\ALPHABET B \subseteq \mathbb R$, and any
      realization $i_t$ of $I_t$ and any choice $u_t$ of $U_t$, we have that
      \begin{equation}
        \label{eq:info-prediction}
        \PR(R_t \in \ALPHABET B \mid I_t = i_t, U_t = u_t)
         =
        \PR(R_t \in \ALPHABET B \mid Z_t = f_t(i_t), U_t = u_t).
      \end{equation}
  \end{enumerate}
\end{definition}
Based on the properties of the information state, we can prove the following.
\begin{theorem}[Structure of optimal control laws]
  \label{thm:structure}
  Let $\{Z_t\}_{t=0}^\infty$, $Z_t \in \ALPHABET Z_t$, be an information-state
  process. Then,
  \begin{enumerate}
    \item The information state absorbs the effect of available information on
      expected future rewards, i.e., for any realization $i_t$ of the
      information state $I_t$, any choice $u_t$ of $U_t$ and any choice of
      future strategy $\VEC g_{(t)} = (g_{t+1}, g_{t+2}, \dots)$, we have that
      \begin{equation}
        \label{eq:prop-1}
        \EXP^{\VEC g_{(t)}}\Big[
          \sum_{\tau=t}^\infty \beta^\tau R_\tau \Bigm| I_t = i_t, U_t = u_t \Big]
          =
          \EXP^{\VEC g_{(t)}}\Big[
            \sum_{\tau=t}^\infty \beta^\tau R_\tau \Bigm| Z_t = f_t(i_t), U_t = u_t \Big].
      \end{equation}
    \item Therefore, $Z_t$ is a sufficient statitistic for performance
      evaluation and there is no loss of optimality in restricting attention to
      control laws of the form $g_t \colon Z_t \mapsto U_t$.
  \end{enumerate}
\end{theorem}

\begin{theorem}[Dynamic programming decomposition]
  \label{thm:DP}
  Assume that the probability distributions in the right-hand side
  of~\eqref{eq:d-dynamics}, \eqref{eq:d-observation}, \eqref{eq:info-evolution}
  and~\eqref{eq:info-prediction} are time-invariant. Let $\{Z_t\}_{t=0}^\infty$,
  be an information-state process such that the space of realization of $Z_t$ is
  time-invariant, i.e., $Z_t \in \ALPHABET Z$.
  \begin{enumerate}
    \item For any choice of future strategy $\VEC g_{(t)} = (g_{t+1}, g_{t+2},
      \dots)$, where $g_\tau$, $\tau > t$,  is of the form $g_\tau \colon Z_\tau
      \mapsto U_\tau$ and for any realization $z_t$ of $Z_t$ and any choice
      $u_t$ of $U_t$, we have that
      \begin{multline}
        \label{eq:prop-2}
        \EXP^{\VEC g_{(t)}}\bigg[
          \EXP^{\VEC g_{(t+1)}}\Big[
            \sum_{\tau=t+1}^\infty \beta^\tau R_\tau \Bigm| Z_{t+1}, U_{t+1} =
            g_{t+1}(Z_{t+1})\Big] \bigg| Z_t = z_t, U_t = u_t \bigg]
        \\
        =
        \EXP^{\VEC g_{(t)}}\Big[
            \sum_{\tau=t+1}^\infty \beta^\tau R_\tau 
            \Big| Z_t = z_t, U_t = u_t \Big]
      \end{multline}
    \item There exists a
      time-invariant optimal strategy $\VEC g^* = (g^*, g^*, \dots)$ that is
      given by
      \begin{subequations}
        \label{eq:DP}
        \begin{equation}
          \label{eq:DPa}
          g^*(z) = \arg \sup_{u \in \ALPHABET U} Q(z,u),
          \quad \forall z \in \ALPHABET Z
        \end{equation}
      where $Q$ is the fixed point solution of the following \emph{dynamic
      program}%
      \footnote{In general, a dynamic program may not have an unique
        solution, or any solution at all. In this paper, we ignore the issue
        of existence of such a solution and refer the reader
        to~\cite{HernandezLermaLasserre:1996} for details.}
        \begin{align}
          Q(z,u) &= \EXP[ R_{t} + \beta V(Z_{t+1}) \mid Z_t = z, U_t = u ],
          \quad \forall z \in \ALPHABET Z, \ u \in \ALPHABET U;
          \\
          V(z)   &= \sup_{u \in \ALPHABET U} Q(z,u),
          \quad \forall z \in \ALPHABET Z.
          \label{eq:DPc}
        \end{align}
      \end{subequations}
  \end{enumerate}
\end{theorem}

The dynamic program can be solved using different methods such as
value-iteration, policy-iteration, or linear-programming.
See~\cite{Puterman:1994} for details.

\begin{remark}
  Identifying an appropriate information-state process for a system resolves the
  two conceptual difficulties described in Sec.~\ref{sec:difficulties}.
  Instead of solving a functional optimization problem to find the optimal
  infinite sequence of control laws, we only need to solve a set of parametric
  optimization problems to find the best control action for each realization of
  information state. A solution to these set of equations determines a control
  law $g^* \colon z \mapsto u$ such that the time-invariant strategy $\VEC g^* =
  (g^*, g^*, \dots)$ is globally optimal. So, we only need to implement one
  control law $g^*$ to implement an optimal control strategy. 
\end{remark}

\begin{remark}
  An important property of the information state is that the conditional future
  cost, which is given by~\eqref{eq:prop-1}, does not depend on the past and
  current control strategy $(g_0, g_1, \dots, g_t)$. This \emph{strategy
  independence} of future cost is critical to obtain a
  recurrence relation for the conditional future cost~\eqref{eq:prop-2} that does not
  depend on the current control law $g_t$. Based on this
  recurrence, we can convert the functional optimization problem of finding the
  best control law $g_t$ into a set of parametric optimization problem of
  finding the best control action $U_t$ for each realization of the information
  state $Z_t$. One of the difficulties in decentralized stochastic control is
  that the future conditional cost from the point of view of a controller
  depends on the \emph{past} and the future control strategy $\VEC g$.
  Therefore, it is not possible to get a dynamic programming decomposition where
  each step is a parametric optimization problem. 
\end{remark}

\begin{remark}
  The information-state based solution approach presented above is equivalent to
  the standard description of centralized stochastic control. In particular:
  \begin{enumerate}
    \item In a \emph{Markov decision process (MDP)}, the controller perfectly
      observes the state process, i.e., $Y_t = X_t$ or equivalently the state
      $X_t$ is a function of the information $I_t$. For such a system $Z_t =
      X_t$ is an information state.

    \item In a \emph{partially observable Markov decision process (POMDP)}, the
      belief state $Z_t(\cdot) = \PR(X_t = \cdot \mid I_t)$ is always an
      information state.
  \end{enumerate}
  In general, a system may have more than one information-state process;
  Theorems~\ref{thm:structure} and~\ref{thm:DP} hold for any information-state
  process. In the next section, we present an example that illustrates, among
  other things, why one information-state process may be preferable to another.
\end{remark}


\subsection {An example}
\label{sec:c-example}

To illustrate the concepts described above, consider an example of a device
transmitting over a communication channel. This may be considered as a special
case of the example of Sec.~\ref{sec:d-example} in which one of the devices never
transmits.

\begin{description}
  \item[\emph{Packet arrival at the device.}]
    Packets arrive at the device according to a Bernoulli process
    $\{W_t\}_{t=0}^\infty$ with rate $p$, i.e., $W_t \in \{0,1\}$ and $\PR(W_t = 1
    \mid W_{1:t-1}) = p$. The device may store $N_t \in \{0,1\}$ packets in a
    buffer. If a packet arrives when a buffer is full, the packet is dropped.

  \item[\emph{Channel model.}]
    The channel model is exactly same as that of Sec.~\ref{sec:d-example}.

  \item[\emph{System dynamics.}]
    At time $t$, the device transmits $U_t \in \{0,1\}$ packets, $U_t \le N_t$.
    If the device transmits when the channel is idle, the transmission is
    successful and the transmitted packet is removed from the buffer. Otherwise,
    the transmission is unsuccessful. Thus, the state of the buffer evolves as
    \[
      N_{t+1} = \min\{ N_t - U_t (1 - S_t) + W_t, 1\}.
    \]
    Each transmission costs $c$ and a successful transmission yields a reward
    $r$. Thus, the total reward is given by
    \[
      R_t = U_t [-c + r (1-S_t)].
    \]

  \item [\emph{Observation model.}]
    The controller perfectly observes the number $N_t$ of packets in the buffer. In
    addition, it observes a channel-state \emph{only if it
    transmits}. Let $H_t$ denote this additional observation. Then $H_t =
    S_{t-1}$ if $U_{t-1} = 1$, otherwise $H_t = \mathfrak E$ (which denotes no
    observation).

  \item[\emph{Information structure.}]
    The information $I_t$ available at the device is given by $I_t = \{N_{1:t},
    U_{1:t-1}, H_{1:t}\}$. The device chooses $U_t$ using a control law $g_t
    \colon I_t \mapsto U_t$. The objective is to pick a control strategy $\VEC g
    = (g_0, g_1, \dots)$ to maximize the expected discounted reward.
\end{description}

The model described above is a centralized stochastic control system with state $X_t =
(N_t, S_t)$, observation $Y_t = (N_t, H_t)$, reward $R_t$, and control $U_t$;
one may verify that these processes satisfy~\eqref{eq:d-dynamics}
and~\eqref{eq:d-observation} (with $n = 1$).

Let $\xi_t \in [0,1]$ denote the posterior probability that the channel is
busy, i.e.,
\[
  \xi_t \DEFINED \PR(S_t = 1 \mid H_{1:t}).
\]
One may verify that $Z_t = (N_t, \xi_t)$ is an information state that
satisfies~\eqref{eq:info-evolution} and~\eqref{eq:info-prediction}. So, there is
no loss of optimality in using control laws of the form $g_t : (N_t, \xi_t)
\mapsto U_t$. The information state takes value in the uncountable space
$\{0,1\} \times [0,1]$. Since $\xi_t$ is a posterior distribution, we can use
the computational techniques of POMDPs~\cite{Zhang:phd,ShaniPineauKaplow:2013}
to numerically solve the corresponding dynamic program. 

However, a simpler dynamic programming decomposition is possible by
characterizing the reachable set of $\xi_t$, which is given by 
\begin{subequations}
  \label{eq:xi-state}
\begin{equation}
  \ALPHABET Q \DEFINED \{ q_{0,k} \mid k \in \naturalnumbers \}
  \cup
  \{ q_{1,k} \mid k \in \naturalnumbers \}
\end{equation}
where
\begin{equation}
  q_{s,k} \DEFINED \PR(S_k = 1 \mid S_0 = s), \quad
  \forall s \in \{0,1\}, \ k \in \naturalnumbers.
\end{equation}
\end{subequations}
Therefore, $\{(N_t, \xi_t)\}_{t=0}^\infty$, $(N_t, \xi_t) \in \{0, 1\} \times
\ALPHABET Q$, is an alternative information-state process. In this alternative
characterization, the information state is denumerable and we may use
finite-state approximations to solve the corresponding dynamic
program~\cite{White:1980,HernandezLerma:1986,CavazosCadena:1986,Flam:1987}.

The dynamic program for this alternative characterization is given below. Let
$\overline p = 1 - p$ and $\overline q_{s,k} = 1 - q_{s,k}$. Then for $s \in \{0,1\}$ and $k \in
\naturalnumbers$, we have that\footnote{Note that $\{ q_{s,k} \mid s \in \{0,1\}$ and $k \in
\naturalnumbers \}$ is equivalent to the reachable set $\ALPHABET Q$ of
$\xi_t$.}
\begin{subequations}
  \label{eq:c-DP-ex}
  \begin{align}
    V(0, q_{s,k}) &= \beta
      \big[ \overline p V(0, q_{s,k+1}) + p V(1, q_{s,k+1}) \big]
      \\
    V(1, q_{s,k}) &= \max\big\{
      \beta V(1, q_{s,k+1}), \notag
      \\
      & \qquad
      \overline q_{s,k} r - c + \beta
      \big[ \overline p \, \overline q_{s,k} V(0, q_{0,1})
        + (p \, \overline q_{s,k} + q_{s,k})  V(1, q_{s,k+1}) \big]
    \big\}
    \label{eq:s-DP-2}
  \end{align}
\end{subequations}
where the first alternative in the right hand side of~\eqref{eq:s-DP-2}
corresponds to choosing $u = 0$ while the second corresponds to choosing $u =
1$.  The resulting optimal strategy for $\beta = 0.9$, $\alpha_0 = \alpha_1 =
0.75$, $r = 1$ and various values of $c$ and $p$ is shown in
Table~\ref{tab:c-structure}.

\begin{table}
  \centering
  \caption{Optimal strategy for the example of Sec.~\ref{sec:c-example} for
$\beta = 0.9$, $\alpha_0 = \alpha_1 = 0.75$, $r=1$ and various values of $c$ and $p$. To
succinctly represent the optimal strategy, each cell shows $(k_0, k_1)$ where
$k_s = \min\{ k \in \naturalnumbers \mid g(1, q_{s,k}) = 1 \}$, for $s \in
\{0,1\}$. }
  \label{tab:c-structure}
  \begin{tabular}{c c c c c c}
    \toprule
    & c = 0.1 & c = 0.2 & c = 0.3 & c = 0.4 & c = 0.5 \\
    \midrule
    p = 0.1 & (1,1) & (1,2) & (1,3) & (1,4) & (1,7) \\
    p = 0.2 & (1,1) & (1,2) & (1,3) & (1,4) & (1,6) \\
    p = 0.3 & (1,1) & (1,2) & (1,2) & (1,3) & (1,5) \\
    p = 0.4 & (1,1) & (1,2) & (1,2) & (1,3) & (1,4) \\
    \bottomrule
  \end{tabular}
\end{table}

As is illustrated by the above example, a general solution methodology for
centralized stochastic control is as follows:
\begin{enumerate}
  \item Identify an information-state process for the given system.
  \item Obtain a dynamic program corresponding to the information-state process.
  \item Either obtain an exact analytic solution of the dynamic program (which
    is only possible for simple stylized models), or obtain an approximate
    numerical solution of the dynamic program (as was done in the example
    above), or prove qualitative properties of the optimal solution (e.g., in
    the above example, for appropriate values of $c$, $r$, and $\VEC P$, the set
    $T(s,n) = \{ k \in \naturalnumbers \mid g^*(n, q_{s,k}) = 1 \}$ is convex).
\end{enumerate}
In the rest of this paper, we explore whether a similar solution approach is
possible for decentralized stochastic control problems.

\section {Conceptual difficulties in dynamic programming for decentralized stochastic control}
\label{sec:discussion}

Recall the two conceptual difficulties that arise in  decentralized stochastic
control and were described in Sec.~\ref{sec:difficulties}. Similar difficulties
arise in centralized stochastic control, where they are resolved by identifying
an appropriate information-state process. It is natural to ask if a similar
simplification is possible in decentralized stochastic control. In particular:
\begin{enumerate}
  \item Is it possible to identify an information state $Z^i_t$, $Z^i_t \in
    \ALPHABET Z^i_t$, such that there is no loss of optimality in restricting
    attention to controllers of the form $g^i_t \colon Z^i_t \mapsto U^i_t$?
  \item If the probability distributions in the right hand side
    of~\eqref{eq:d-dynamics} and~\eqref{eq:d-observation} are time-invariant,
    is it possible to identify a dynamic programming decomposition that
    determines optimal control strategies for all controllers?
\end{enumerate}
The second question is significantly more important, and considerably harder,
than the first. There are two approaches to find a dynamic programming
decomposition. The first approach is to find a set of coupled dynamic programs,
where each dynamic program is associated with a controller and determines the
``optimal'' control strategy at that controller. The second approach is to find
a dynamic program that simultaneously determines the optimal control strategy at
all controllers. 

It is not obvious how to identify such dynamic programs. Let's conduct a thought
experiment in which we assume that such dynamic programs have been identified
and let's try to identify the implications. The description below is qualitative;
the mathematical justification is presented later in the paper. 

Consider the first approach.
Suppose we are able to find a set of coupled dynamic programs, where the dynamic
program for controller~$i$, which we refer to as $\mathcal D^i$, determines the
``optimal'' strategy $\VEC g^i$ for controller~$i$. We use the term optimal in
quotes because we cannot isolate an optimization problem for controller~$i$
until we specify the control strategy $\VEC g^{-i}$ for all other controllers.
Therefore, dynamic program $\mathcal D^i$ determines the \emph{best response
strategy} $\VEC g^i$ for a particular choice of control strategies $\VEC g^{-i}$
for other controllers. With a slight abuse of notation, we can write this as
\[
  \VEC g^i = \mathcal D^i(\VEC g^{-i}).
\]
Any solution $\VEC g^{*} = (\VEC g^{*,1}, \dots, \VEC g^{*,n})$ of these coupled
dynamic programs will have the property that for any controller~$i$, $i \in
\{1,\dots,n\}$, given that all other controllers are using the strategy $\VEC
g^{*,-i}$, controller~$i$ is playing its best response strategy $\VEC g^{*,i} =
\mathcal D^i(\VEC g^{*,-i})$. Such a strategy is called a \emph{person-by-person
optimal} strategy (which is related to the notion of local optimum in
optimization theory and the notion of Nash equilibrium in game theory). In general,
a person-by-person optimal strategy need not be globally optimal; in fact, a
person-by-person strategy may perform arbitrarily bad as compared to the
globally optimal strategy. In conclusion, unless we impose further restrictions
on the model, a set of coupled dynamic programs cannot determine a globally
optimal strategy. 

Now, consider the second approach. Suppose we are able to find a dynamic
program similar to~\eqref{eq:DPa}--\eqref{eq:DPc} that determines the optimal
control strategies for all controllers. All controllers must be able to use this
dynamic program to find their control strategy. Therefore, the information-state
process $\{Z_t\}_{t=0}^\infty$ of such a dynamic program must have the following
property: \emph{$Z_t$ is a function of the information $I^i_t$ available to
every controller $i$, $i \in \{1,\dots,n\}$}. In other words, the information
state must be measurable with respect to the \emph{common knowledge} (in the
sense of Aumann~\cite{Aumann:1976}) between the controllers.

In centralized stochastic control, we first showed that there was no loss of
optimality in restricting attention to control laws of the form $g \colon Z
\mapsto U$, then used this in step~\eqref{eq:DPc} to convert the functional
optimization problem of finding the best control law $g$ into a parametric
optimization problem of finding the best $u$ for each realization $z$ of
information-state. In the decentralized case, we just argued that the
information state $Z_t$ must be commonly known to all controllers. Therefore, if
we restrict attention to control laws of the form $g^i_t \colon Z_t \mapsto
U^i_t$, then each controller would be ignoring its \emph{local information}
(i.e., the information not commonly known to all controllers). Hence, a
restriction to control laws of the form $g^i_t \colon Z_t \mapsto U^i_t$ cannot
be without loss of optimality. 

If we have a dynamic program similar to~\eqref{eq:DPa}--\eqref{eq:DPc} that
uses information-state process $\{Z_t\}_{t=0}^\infty$ to determine the optimal
control strategy for all controllers, then restricting attention to control
laws of the form $g^i_t \colon Z_t \mapsto U^i_t$ will result in loss of optimality.
Therefore, the step corresponding to~\eqref{eq:DPc} cannot be a parametric
optimization problem and it must be a functional optimization problem. 

Now let's try to characterize the nature of the functional optimization problem
corresponding to~\eqref{eq:DPc}. The only way in which the solution of this
functional optimization problem will determine optimal control strategies for
all controllers is as follows. Let $L^i_t$ denote the \emph{local information}
at each controller so that $Z_t$ and $L^i_t$ are sufficient to determine
$I^i_t$. Then, for a particular realization $z$ of the information-state, the
step corresponding to~\eqref{eq:DPc} of the dynamic program must determine
functions $(\gamma^1, \dots, \gamma^n)$ such that: (i)~$\gamma^i$ gives
instructions to controller~$i$ on how to use its local information $L^i_t$ to determine the
control action $U^i_t$; and (ii)~computing $(\gamma^1, \dots, \gamma^n)$ for
each realization of the information state is equivalent to choosing $(g^1,
\dots, g^n)$. These steps are made precise in Sec.~\ref{sec:common}. 

The above discussion shows that dynamic programming for decentralized stochastic
control will be different from that for centralized stochastic control. Either
we must be content with a person-by-person optimal strategy; or, if we pursue
global optimality, then we must be willing to solve functional optimization
problems in the step corresponding to~\eqref{eq:DPc} in an appropriate dynamic program.

In the literature, the first approach is called the \emph{person-by-person
approach} and the second approach is called the \emph{common-information
approach}. We describe both these approaches in the next section.

\section {The person-by-person approach}
\label{sec:p-by-p}

The person-by-person approach is motivated by the computational approaches for
finding Nash equilibrium in game theory. It was proposed by Marschak and
Radner~\cite{Radner:1962,MarschakRadner:1972} in the context of static systems
with multiple controllers and has been subsequently used in dynamic systems as
well. The main idea behind the person-by-person approach is to decompose the
decentralized stochastic control problem into a series of centralized stochastic
control sub-problems, each from the point of view of a single controller, and
use the solution techniques of centralized stochastic control to simplify these
sub-problems. The person-by-person approach is used to identify structural
results as well as identify coupled dynamic programs to find person-by-person
optimal (or equilibrium) strategies. 

\subsection {Structure of optimal control strategies}

To find the structural results, proceed as follows. Pick a controller that has
perfect recall, say~$i$;
\emph{arbitrarily} fix the control strategies $\VEC g^{-i}$ of all controllers
except controller~$i$ and consider the sub-problem of finding the \emph{best
response} strategy $\VEC g^i$ at controller~$i$. Since controller~$i$ has
perfect recall, this sub-problem is centralized. Suppose that we identify an
information-state process $\{\tilde I^i_t\}_{t=0}^\infty$ for this sub-problem.
Then, there is no loss of (best-response) optimality in restricting attention to control laws of
the form $\tilde g^i_t \colon \tilde I^i_t \to U^i_t$ at controller~$i$.

Recall that the choice of control strategies $\VEC g^{-i}$ was completely
arbitrary. Suppose the structure of $\tilde g^i_t$ does not depend on the choice
of control strategies $\VEC g^{-i}$ of other controllers, then there is no loss
of (global) optimality in restricting attention to control laws of the form
$\tilde g^i_t$ at controller~$i$. 

Repeat this procedure at all controllers that have perfect recall. Let $\{\tilde I^i_t\}_{t=0}^\infty$
be the information-state processes identified at controller~$i$, $i \in
\{1,\dots,n\}$. Then there is no loss of global optimality in restricting
attention to the information structure $(\tilde I^i_t, i \in \{1,\dots,n\},
t=0,1,\dots)$. 

To illustrate this approach, consider the example of the decentralized control
system of Sec.~\ref{sec:d-example}. Arbitrarily fix the control strategy $\VEC
g^j$ of controller~$j$, $j \in \{1,2\}$, and consider the sub-problem of finding
the best response strategy $\VEC g^i$ of controller~$i$, $i = 3 - j$. Since
controller~$i$ has perfect recall, the subproblem of finding the best response
strategy $\VEC g^i$ is a centralized stochastic control problem. To simplify
this centralized stochastic control problem, we need to identify an information
state as described in Definition~\ref{def:info}

Recall assumption (A) that states that the packet-arrival process at the two
devices are independent. Under this assumption, we can show that
\begin{multline}
  \label{eq:d-indep}
  \PR(N^1_{1:t}, N^2_{1:t} \mid H_{1:t}, U^1_{1:t-1}, U^2_{1:t-1})
  \\=
  \PR(N^1_{1:t} \mid H_{1:t}, U^1_{1:t-1}, U^2_{1:t-1})
  \PR(N^2_{1:t} \mid H_{1:t}, U^1_{1:t-1}, U^2_{1:t-1})
\end{multline}
Using the above conditional independence, we can show that for any choice of
control strategy $\VEC g^j$, $\tilde I^i_t = \{N^i_t, H_{1:t}, U^1_{1:t-1},
U^2_{1:t-1}\}$ is an information state for controller~$i$. By
Theorem~\ref{thm:structure}, we get that there is no
loss of optimality (for best response strategy) in restricting attention to
control laws of the form $\tilde g^i_t \colon \tilde I^i_t \mapsto U^i_t$. Since the
structure of the optimal best response strategy does not depend on the choice of
$\VEC g^j$, there is no loss of global optimality in restricting attention to
control laws of the form $\tilde g^i_t$. Equivalently, there is no loss of
optimality in assuming that the system has a simplified information structure
$(\tilde I^i_t, i \in \{1,2\}, t = 0, 1, \dots)$. 

\subsection {Coupled dynamic program for person-by-person optimal solution}

Based on the discussion in Sec.~\ref{sec:discussion}, it is natural to ask if
the method described above can be extended to find coupled dynamic programs that
determine person-by-person optimal strategies when all controllers have perfect
recall and the model is time-homogeneous, i.e., the probability distribution on
the right hand side of~\eqref{eq:d-dynamics} and~\eqref{eq:d-observation} are
time-invariant.

Suppose that by using the person-by-person approach, we find that there is no loss
of optimality in restricting attention to the information structure $(\tilde
I^i_t, i \in \{1,\dots, n\}, t = 0,1,\dots)$ and control strategies $\tilde
g^i_t \colon \tilde I^i_t \mapsto U^i_t$, $i \in \{1, \dots, n\}$. Pick a
controller, say~$i$, and arbitrarily fix the control strategies $\tilde {\VEC
g}^{-i}$ of all controllers other than~$i$. Is it possible to use
Theorem~\ref{thm:DP} to find the best response strategy $\tilde {\VEC g}^i$ at
controller~$i$? In general, the answer is no because of the following reasons.
\begin{enumerate}
  \item The information-state process $\{\tilde I^i_t\}_{t=0}^\infty$ in general
    does not take values in time-invariant space (e.g., in the above example,
    $\tilde I_t = \{N^i_t, H_{1:t}, U^1_{1:t-1}, U^2_{1:t-1}\}$). A fortiori, we
    cannot show that restricting attention to time-invariant strategies is
    without loss of optimality.

  \item Assume that for every controller~$i$, $i \in \{1,\dots,n\}$, the
    information-state process $\{\tilde I^i_t\}_{t=0}^\infty$, takes value in a
    time-invariant space. Even then, when we \emph{arbitrarily} fix the control
    strategies $\tilde {\VEC g}^{-i}$ of all other controllers, the dynamical
    model seen by controller~$i$ is not time homogeneous. For the dynamic
    model from the point of view of controller~$i$ to be time-homogeneous, we
    must further assume that each controller~$j$, $j \neq i$, is using a
    time-invariant strategy $\tilde {\VEC g}^j$. 
\end{enumerate}

Therefore, if the information-state process $\{\tilde I^i_t\}_{t=0}^\infty$
for every controller~$i$ takes value in time-invariant space and we \emph{a
priori} restrict attention to time-invariant strategies (even if such a
restriction results in loss of optimality), then the problem of finding the
best-response strategy at a particular controller is a \emph{time-homogeneous}
expected discounted cost problem that can be solved using Theorem~\ref{thm:DP}. 

In particular, let $\ALPHABET D^i$ denote the dynamic program to find the best
response strategy $\tilde {\VEC g}^i$ for controller~$i$ when all other
controllers are using a \emph{time-invariant} strategy $\tilde {\VEC g}^j =
(\tilde g^j, \tilde g^j, \dots)$, $j \neq i$. From Theorem~\ref{thm:DP} we know
that $\tilde {\VEC g}^i$ is also time-invariant. With a slight abuse of
notation, we denote this relationship as
\[
  \tilde g^i = \mathcal D^i(\tilde g^j, j \neq i).
\]
We can write similar dynamic programs for all controllers~$i$, $i = \{1, \dots,
n\}$, giving $n$ coupled dynamic programs. 

As described in Sec.~\ref{sec:discussion}, a solution $(g^{*,i}, i \in
\{1,\dots,n\})$ of these coupled dynamic programs is a \emph{person-by-person}
optimal strategy. Such a time-invariant person-by-person optimal strategy need
not be globally optimal for two reasons. Firstly, there might be other
time-invariant person-by-person strategies that achieve a higher expected
discounted reward. Secondly, we haven't showed that restricting attention to
time-invariant strategies is without loss of optimality. Thus, there might be
other time-\emph{varying} strategies that achieve higher expected discounted
reward. 

Such coupled dynamic programs have been used to find person-by-person optimal
strategies in sequential detection
problems~\cite{TeneketzisHo:1987,TeneketzisVaraiya:1984}.

\section {The common-information approach}
\label{sec:common}

The common-information approach was proposed by Nayyar, Mahajan,
Teneketzis~\cite{Nayyar:phd,MNT:tractable-allerton,NMT:partial-history-sharing,NMT:CI-chapter}
and provides a dynamic programming decomposition for for a subclass of decentralized control
systems that determines optimal control strategies for all controllers.
Variation of this approach had been used for specific information structures
including delayed state sharing~\cite{AicardiDavoliMinciardi:1987}, partially
nested systems with common
past~\cite{CasalinoDavoliMinciardiPuliafitoZoppoli:1984}, teams with sequential
partitions~\cite{Yoshikawa:1978}, periodic sharing information
structure~\cite{OoiVerboutLudwigWornell:1997}, and belief sharing information
structure~\cite{Yuksel:2009}.

This approach formalizes the intuition presented in Sec.~\ref{sec:discussion}:
to obtain a dynamic program that determines optimal control strategies for all
controllers, the information-process must be measurable at all controllers and,
at each step of the dynamic program, we must solve a functional optimization
problem that determines instructions to map local information to control action
for each realization of the information state. 

To formally describe this intuition, split the information available at each
controller into two parts: the \emph{common information}
\[
  C_t = \bigcap_{\tau \ge t} \bigcap_{i=1}^n I^i_\tau
\]
and the \emph{local information}
\[
  L^i_t = I^i_t \setminus C_t,
  \quad \forall i \in \{1, \dots, n\}.
\]

By construction, the common and local information determine the total
information, i.e., $I^i_t = C_t \cup L^i_t$ and the common information is
nested, i.e., $C_t \subseteq C_{t+1}$. 

For simplicity of presentation, we restrict to \emph{partial history sharing}
information structures~\cite{NMT:partial-history-sharing,NMT:CI-chapter}. The
common information approach is applicable to a more general class of models.
See~\cite{Nayyar:phd} for details.

\begin{definition} \label{def:PHS}
  An information structure is called \emph{partial history sharing} information
  structure when the following conditions are satisfied:
  \begin{enumerate}
    \item For any set of realizations $\ALPHABET A$ of $L^i_{t+1}$ and any
      realization $c_t$ of $C_t$, $\ell^i_t$ of $L^i_t$, $u^i_t$ of $U^i_t$ and
      $y^i_{t+1}$ of $Y^i_{t+1}$, we have
      \begin{multline*}
        \PR(L^i_{t+1} \in \ALPHABET A \mid
            C_t = c_t, L^i_t = \ell^i_t, U^i_t = u^i_t, Y^i_{t+1} = y^i_{t+1})
        \\ = 
        \PR(L^i_{t+1} \in \ALPHABET A \mid
            L^i_t = \ell^i_t, U^i_t = u^i_t, Y^i_{t+1} = y^i_{t+1})
      \end{multline*}
    \item the size of the local information is uniformly bounded\footnote{This
        condition is needed to ensure that the information-state is
      time-invariant and, as such, may be ignored for finite horizon
      models~\cite{NMT:partial-history-sharing}.}, i.e., there
      exists a $k$ such that for all $t$ and all $i \in \{1, \dots, n\}$,
      $|\ALPHABET L^i_t| \le k$, where $\ALPHABET L^i_t$ denotes the space of
      realizations of $L^i_t$. 
  \end{enumerate}
\end{definition}

Systems with complete information sharing, $k$-step delayed sharing, and
$k$-step periodic sharing information structures described in
Sec.~\ref{sec:info-structure} are special cases of partial history sharing information
structures. The model of Sec.~\ref{sec:d-example} does not have a partial
history sharing structure, but when we restrict attention to the information
structure $(\tilde I^i_t, i \in \{1,2\}, t = 0,1,\dots)$ where $\tilde I^i_t =
\{N^i_t, H_{1:t}, U^1_{1:t-1}, U^2_{1:t-1}\}$, then the model has partial
history sharing information structure. 

The objective of the common-information approach is to identify a dynamic
program that determines optimal control strategies for all controllers. The
simplest way to describe the approach is to construct a centralized stochastic
control problem that gives rise to such a dynamic program. We can convert a
decentralized stochastic control problem into a centralized stochastic control
problem by exploiting the fact that planning is centralized, i.e., the control
strategies for for all controllers are chosen before the system starts running and,
therefore, optimal strategies can be searched in a centralized manner.

The construction of an appropriate dynamic program relies on partial evaluation
of a function defined below.
\begin{definition}
  For any function $f \colon (x,y) \mapsto z$ and a value $x_0$ of $x$,
  the \emph{partial evaluation} of $f$ and $x = x_0$ is a function $g \colon y
  \mapsto z$ such that for all values of $y$,
  \[
    g(y) = f(x_0, y).
  \]
\end{definition}
For example, if $f(x,y) = x^2 + xy + y^2$, then the partial evaluation of $f$ at
$x = 2$ is $g(y) = y^2 + 2y + 4$.

\newlength\normalparindent
\setlength\normalparindent{15pt}

The common-information approach proceeds as
follows~\cite{NMT:partial-history-sharing, NMT:CI-chapter}:
\begin{enumerate}
  \item \emph{Construct a centralized coordinated system}.

    The first step of the common-information approach is to construct a
    centralized stochastic control system called the \emph{coordinated system}.
    The controller of this system, called the \emph{coordinator}, observes the
    common information $C_t$ and chooses the partially evaluated control laws
    $g^i_t$, $i \in \{1,\dots, n\}$ at $C_t$. Denote these partial evaluations
    by $\Gamma^i_t$ and call them \emph{prescriptions}. These prescriptions tell
    the controllers how to map their local information information into control
    actions; in particular $U^i_t = \Gamma^i_t(L^i_t)$. The decision rule
    $\psi_t \colon C_t \mapsto (\Gamma^1_t, \dots, \Gamma^n_t)$ to choose the
    prescriptions is called a \emph{coordination law}. 
    
    \hspace* \normalparindent
    The coordinated system has only one controller, the coordinator, which has
    perfect recall; the controllers of the original system are passive agents
    that simply use the prescriptions given by the coordinator. Hence, the
    coordinated system is a centralized stochastic control system with the state
    process $\{(X_t, L^1_t, \dots, L^n_t)\}_{t=0}^\infty$, the observation
    process $\{C_t\}_{t=0}^\infty$, the reward process $\{R_t\}_{t=0}^\infty$,
    and the control process $\{(\Gamma^1_t, \dots, \Gamma^n_t)\}_{t=0}^\infty$.

    \hspace* \normalparindent
    In contrast to centralized stochastic control, the control process of the
    coordinated system is a sequence of functions. Consequently, when we
    describe the dynamic program to find the best ``control action'' for each
    realization of the information state, the step corresponding
    to~\eqref{eq:DPc} will be a functional optimization problem. 

  \item \emph{Simplify the coordinated system}

    Let $\{Z_t\}_{t=0}^\infty$, $Z_t \in \ALPHABET Z_t$, be an information-state
    process for the coordinated system.\footnote{Since the coordinated system is
      a POMDP, the process $\{\pi_t\}_{t=0}^\infty$, where $\pi_t$ is the
      conditional probability measure on $(X_t, L^1_t, \dots, L^n_t)$
    conditioned on $C_t$, is always an information-state process.} Therefore,
    there is no loss of optimality in restricting attention to coordination
    laws of the form
    \[
      \psi_t \colon Z_t \mapsto (\Gamma^1_t, \dots, \Gamma^n_t).
    \]
    When the probability distributions on the right hand side
    of~\eqref{eq:d-dynamics} and~\eqref{eq:d-observation} are time-invariant,
    the evolution of $Z_t$ is time-invariant, 
    and the state space $\ALPHABET Z_t$ of the realizations of $Z_t$ is
    time-invariant, i.e., $\ALPHABET Z_t = \ALPHABET Z$,
    then there exists a time-invariant coordination strategy $\boldsymbol
    \psi^* = (\psi^*, \psi^*, \dots)$ where $\psi^*$ is given by 
    \begin{subequations}
    \begin{equation}
      \label{eq:d-DPa}
      \psi^*(z) = \arg \sup_{(\gamma^1,\dots,\gamma^n)} 
        Q(z, (\gamma^1,\dots,\gamma^n)), \quad \forall z \in \ALPHABET Z
    \end{equation}
    where $Q$ is the unique fixed point of the following set of equations
      \begin{align}
        Q(z, (\gamma^1,\dots,\gamma^n)) &= \EXP[R_t + \beta V(Z_{t+1}) 
          \mid Z_t = z, \notag \\
          & \qquad\quad 
          \Gamma^1_t = \gamma^1, \dots, \Gamma^n_t = \gamma^n ],
          \quad \forall z \in \ALPHABET Z,\ \forall (\gamma^1,\dots,\gamma^n)
          \\
        V(z) &= \sup_{(\gamma^1,\dots,\gamma^n)} Q(z, (\gamma^1,\dots,\gamma^n))
        \label{eq:d-DPc}
      \end{align}
    \end{subequations}

    Step~\eqref{eq:d-DPc} of the above dynamic program is a functional
    optimization problem. In contrast, step~\eqref{eq:DPc} of the dynamic
    program for centralized stochastic control was a parametric optimization
    problem.

  \item \emph{Show equivalence between the original system and the coordinated
      system and translate the solution of the coordinated system to the
    original system}.

    It can be shown that the coordinated system is equivalent to the original
    system~\cite{NMT:partial-history-sharing}. In particular, for any choice of
    the coordination strategy in the coordinated system, there exists a control
    strategy in the original decentralized system that gives the same expected
    discounted reward, and vice-versa. Using this equivalence, we can translate
    the results of the previous step to the original decentralized system.
    Hence, if $\{Z_t\}_{t=0}^\infty$ is an information-state process for the
    coordinated system, then there is no loss of optimality in restricting
    attention to control strategies of the form
    \[
      g^i_t \colon (Z_t, L^i_t) \mapsto U^i_t.
    \]
    Furthermore, if  $\boldsymbol \psi^* = (\psi^*, \psi^*, \dots)$ is an
    optimal time-invariant coordination strategy for the coordinated system, then the
    time-invariant control strategies $\VEC g^{i,*} = (g^{i,*}, g^{i,*},
    \dots)$, $i \in \{1, \dots, n\}$, where
    \[
      g^{i,*}(z, \ell^{i}) = \psi^{i,*}(z)(\ell^{i})
    \]
    and $\psi^{i,*}_t$ denotes the $i$-th component of $\psi_t$, are optimal for
    the original decentralized system.
\end{enumerate}

\begin{remark}
  The coordinated system and the coordinator described above are fictitious and
  used only as a tool to explain the approach. The computations carried out at the
  coordinator are based on the information known to all controllers. Hence, each
  controller can carry out the computations attributed to the coordinator.
  As a consequence, it is possible to describe the above approach without
  considering a coordinator, but in our opinion thinking in terms of a
  fictitious coordinator makes it easier to understand the approach.
\end{remark}

To illustrate this approach, consider the decentralized control example of
Sec.~\ref{sec:d-example}. Start with the simplified information structure
$\tilde I^i_t = \{N^i_t, H_{1:t}, U^1_{1:t-1}, U^2_{1:t-1}\}$ obtained using the
person-by-person approach. The common information is given by 
\[
  C_t = \bigcap_{\tau \ge t} (\tilde I^1_{\tau} \cap \tilde I^2_{\tau}) 
  = \{H_{1:t}, U^1_{1:t-1}, U^2_{1:t-1}\}
\]
and the local information is given by
\[
  L^i_t = \tilde I^i_t \setminus C_t = N^i_t, \quad \forall i \in \{1,2\}.
\]
Thus, in the coordinated system, the coordinator observes $C_t$ and uses the
coordination law $\psi_t \colon C_t \mapsto (\gamma^1_t, \gamma^2_t)$, where
$\gamma^i_t$ maps the local information $N^i_t$ to $U^i_t$. Note that
$\gamma^i_t$ is completely specified by $D^i_t = \gamma^i_t(1)$ because the
constraint $U^i_t \le N^i_t$ implies that $\gamma^i_t(0) = 0$. Therefore, we
may assume that the coordinator uses a coordination law $\psi_t \colon C_t
\mapsto (D^1_t, D^2_t)$, $D^i_t \in \{0,1\}$, $i \in \{1,2\}$ and each device
then chooses a control action according to $U^i_t = N^i_t D^i_t$. The system
dynamics and the reward process are same as in the original decentralized
system.

Since the coordinator has perfect recall, the problem of finding the best
coordination strategy is a centralized stochastic control problem. To simplify
this centralized stochastic control problem, we need to identify an information
state as described in Definition~\ref{def:info}.

Let $\zeta^i_t \in [0,1]$ denote the posterior probability that device~$i$, $i \in
\{1,2\}$ has a packet in its buffer given the channel feedback, i.e.,
\[
  \zeta^i_t = \PR(N^i_t = 1 \mid H_{1:t}, U^1_{1:t-1}, U^2_{1:t-1}), 
  \quad \forall i \in \{1,2\}.
\]
Moreover, as in the centralized case, let $\xi_t \in [0,1]$ denote the posterior
probability that the channel is busy given the channel feedback, i.e., 
\[
  \xi_t = \PR(S_t = 1 \mid H_{1:t}, U^1_{1:t-1}, U^2_{1:t-1}) = 
  \PR(S_t = 1 \mid H_{1:t}).
\]
One may verify that $(\zeta^1_t, \zeta^2_t, \xi_t)$ is an information state that
satisfies~\eqref{eq:info-evolution} and~\eqref{eq:info-prediction}. So, there is
no loss of optimality in using coordination laws of the form $\gamma \colon
(\zeta^1_t, \zeta^2_t, \xi_t) \mapsto (D^1_t, D^2_t)$. This information state
takes values in the uncountable space $[0,1]^3$. Since each component
$\zeta^1_t$, $\zeta^2_t$, and $\xi_t$ of the information state is a posterior
distribution, we can use the computational techniques of
POMDPs~\cite{Zhang:phd,ShaniPineauKaplow:2013} to numerically solve the
corresponding dynamic program.

However, a simpler dynamic programming decomposition is possible by
characterizing the reachable set of the information state. The reachable set of
$\zeta^i_t$ is given by
\begin{subequations}
  \begin{equation}
    \ALPHABET R^i \DEFINED \{ z^i_k \mid k \in \naturalnumbers \} \cup \{ 1 \}
  \end{equation}
  where
  \begin{equation}
    z^i_k \DEFINED \PR(N^i_k = 1 \mid N^i_0 = 0, D^i_{1:k-1} = (0,\dots,0)),
    \quad
    \forall s \in \{0,1\}, \ k \in \naturalnumbers
  \end{equation}
\end{subequations}
and the reachable set of $\xi_t$ is given by $\ALPHABET Q$ defined
in~\eqref{eq:xi-state}. For ease of notation, define $z^i_{\infty} = 1$. 

Therefore, $\{(\zeta^1_t, \zeta^2_t, \xi_t)\}_{t=0}^\infty$, $(\zeta^1_t,
\zeta^2_t, \xi_t) \in \ALPHABET R^1 \times \ALPHABET R^2 \times \ALPHABET Q$, is
an alternative information-state process. In this alternative characterization,
the information state is denumerable and we may use finite-state approximations
to solve the corresponding dynamic
program~\cite{White:1980,HernandezLerma:1986,CavazosCadena:1986,Flam:1987}.

The dynamic program for this alternative characterization is given below. 
Let $\overline q^s_m = 1 - q^s_m$ and $\overline z^i_k = 1 - z^i_k$. Then for $s
\in \{0,1\}$ and $k,\ell \in \naturalnumbers \cup \{\infty\}$ and $m \in
\naturalnumbers$, we have that
\begin{subequations}
  \label{eq:d-DP-ex}
\begin{multline}
  V(z^1_k, z^2_\ell, q^s_m) =
  \max \big\{
    Q(z^1_k, z^2_\ell, q^s_m; (0,0)),
    Q(z^1_k, z^2_\ell, q^s_m; (1,0)),
    \\
    Q(z^1_k, z^2_\ell, q^s_m; (0,1)),
    Q(z^1_k, z^2_\ell, q^s_m; (1,1))
  \big\}
\end{multline}
where $Q(z^1_k, z^2_\ell, q^s_m, (d^1,d^2))$ corresponds to choosing the
prescription $(d^1,d^2)$ and is given by
\begin{align}
    Q(z^1_k, z^2_\ell, q^s_m &; (0,0)) = 
    \beta V(z^1_{k+1}, z^2_{\ell+1}, q^s_{m+1}); 
    \\
    Q(z^1_k, z^2_\ell, q^s_m &; (1,0)) = 
        z^1_k\, \overline q^s_m\, r - z^1_k \, c
     + \beta \Big[
       \overline z^1_k V(z^1_1, z^2_{\ell + 1}, q^s_{m + 1})
      \notag \\
      & \quad\qquad\qquad
      + z^1_k \, \overline q^s_m V( z^1_1, z^2_{\ell + 1}, q^0_1)
      + z^1_k \, q^s_m V( z^1_\infty, z^2_{\ell + 1}, q^1_1)
    \Big]; 
    \\
    Q(z^1_k, z^2_\ell, q^s_m &; (0,1)) = 
        z^2_\ell \, \overline q^s_m\, r - z^2_\ell \, c
     + \beta \Big[
       \overline z^2_\ell V(z^1_{k+1}, z^2_1, q^s_{m + 1})
      \notag \\
      & \quad\qquad\qquad
      + z^2_\ell \, \overline q^s_m V( z^1_{k+1}, z^2_1, q^0_1)
      + z^2_\ell \, q^s_m V( z^1_{k+1}, z^2_\infty, q^1_1)
    \Big];
    \\
    Q(z^1_k, z^2_\ell, q^s_m &; (1,1)) = 
    [z^1_k\, \overline z^2_\ell + \overline z^1_k \, z^2_\ell ] \, \overline q^s_m \, r 
    - [ z^1_k + z^2_\ell ] \, c
    + \beta \Big[
      \overline z^1_k \, \overline z^2_\ell V(z^1_1, z^2_1, q^s_{m+1})
    \notag \\
    & + [ z^1_k \, \overline z^2_\ell + \overline z^1_k \, z^2_\ell ]
      \, \overline q^s_m V(z^1_1, z^2_1, q^0_1)
     + z^1_k z^2_\ell \overline q^s_m V(z^1_\infty, z^2_\infty, q^0_1)
    \notag \\
    & + [ z^1_k + z^2_\ell -  z^1_k \, z^2_\ell ]
      \, q^s_m V(z^1_\infty, z^2_\infty, q^1_1)
      \Big].
\end{align}
\end{subequations}

To describe the optimal strategy, define functions $d$ and $\bar d$ as follows:
\begin{align*}
  d(z^1_k, z^2_\ell) = \begin{cases}
    (1,0), & \text{if $k > \ell$} \\
    (0,1), & \text{if $k < \ell$} \\
    (1,0) \text { or } (0,1), & \text {if $k = \ell$}
  \end{cases}
  \quad \text{and} \quad
  \overline d(z^1_k, z^2_\ell) = \begin{cases}
    (0,1), & \text{if $k > \ell$} \\
    (1,0), & \text{if $k < \ell$} \\
    (1,0) \text { or } (0,1), & \text {if $k = \ell$}
  \end{cases}
\end{align*}
In addition define the sets $\ALPHABET S_n, \hat {\ALPHABET S}_n \subseteq \ALPHABET
R^1 \times \ALPHABET R^2$ for $n \in \naturalnumbers \cup \{\infty\}$ as follows:
\begin{align*}
  S_n &= \{ (z^1_k, z^2_1)    : z^1_k \in \ALPHABET R^1 \text{ and } k \le n \} \cup
        \{ (z^1_1, z^2_\ell) : z^2_k \in \ALPHABET R^2 \text{ and } \ell \le n \}.
  \\
  \hat S_n &= \{ (z^1_k, z^2_\ell) \in \ALPHABET R^1 \times \ALPHABET R^2 :
  \max(k, \ell) \le n \}.
\end{align*}
Using these definitions, define the following functions for $n \in
\naturalnumbers \cup \{\infty\}$.
\begin{enumerate}
  \item $h_n(z^1_k, z^2_\ell) = \begin{cases}
      (1,1), &\text{if $(z^1_k, z^2_\ell) \in S_n$} \\
      d(z^1_k, z^2_\ell), & \text{ otherwise}.
    \end{cases}$

  \item $\hat h_n(z^1_k, z^2_\ell) = \begin{cases}
      (0,0), &\text{if $(z^1_k, z^2_\ell) \in \hat S_n$} \\
      d(z^1_k, z^2_\ell), & \text{ otherwise}.
    \end{cases}$

\end{enumerate}


The optimal strategies obtained by solving~\eqref{eq:d-DP-ex} for $\beta = 0.9$,
$\alpha_0 = \alpha_1 = 0.75$, $r = 1$, $p_1 = p_2 = 0.3$, and different values
of $c$ are given below. 

\begin{enumerate}
  \item When $c = 0.1$ the optimal strategy is given by
    \[
      g^*(z^1_k, z^2_\ell, q^s_m) = \begin{cases}
        h_1(z^1_k, z^2_\ell), & \text{if $s = 0$ and $m = 1$} \\
        h_5(z^1_k, z^2_\ell), & \text{if $s = 1$ and $m = 1$} \\
        h_2(z^1_k, z^2_\ell), & \text{otherwise}.
      \end{cases}
    \]
  \item When $c = 0.2$ the optimal strategy is given by
    \[
      g^*(z^1_k, z^2_\ell, q^s_m) = \begin{cases}
        \overline d(z^1_k, z^2_\ell), & \text{if $s = 1$ and $m = 1$} \\
        d(z^1_k, z^2_\ell), & \text{otherwise}.
      \end{cases}
    \]
  \item  When $c = 0.3$, the optimal strategy is given by
    \[
      g^*(z^1_k, z^2_\ell, q^s_m) = \begin{cases}
        (0,0),              & \text{if $s = 1$ and $m = 1$} \\
        d(z^1_k, z^2_\ell), & \text{otherwise}.
      \end{cases}
    \]
  \item  When $c = 0.4$, the optimal strategy is given by
    \[
      g^*(z^1_k, z^2_\ell, q^s_m) = \begin{cases}
        (0,0),              & \text{if $s = 1$ and $m \le 2$} \\
        d(z^1_k, z^2_\ell), & \text{otherwise}.
      \end{cases}
    \]
  \item  When $c = 0.5$, the optimal strategy is given by
    \[
      g^*(z^1_k, z^2_\ell, q^s_m) = \begin{cases}
        (0,0), & \text{if $s = 1$ and $m \le 3$} \\
        \hat h_1(z^1_k, z^2_\ell), & \text{if $s = 1$, $m = 4$}, \\
        \overline d(z^1_k, z^2_\ell), & \text{if $s = 1$, $m = 5$}, \\
        d(z^1_k, z^2_\ell), & \text{otherwise}.
      \end{cases}
    \]

\end{enumerate}

\begin{remark}
  As we argued in Sec.~\ref{sec:discussion}, if a single dynamic program
  determines the optimal control strategies at all controllers, then the
  step~\eqref{eq:d-DPc} must be a functional optimization problem. Consequently,
  the dynamic program for decentralized stochastic control is significantly
  more difficult to solve than dynamic programs for centralized stochastic
  control. When the observation and control processes are finite valued (as in
  the above example), the space of functions from $L^i_t$ to $U^i_t$ are finite
  and step~\eqref{eq:d-DPc} can be solved by exhaustively searching over all
  alternatives.  
\end{remark}

\begin{remark}
  As in centralized stochastic control, the information-state in decentralized
  control is sensitive to the modeling assumptions. For example, in the above
  example, if we remove assumption~(A), then the conditional independence
  in~\eqref{eq:d-indep} is not valid; therefore, we cannot use the
  person-by-person approach to show that $\{N^i_t, U^1_{1:t-1}, U^2_{1:t-1},
  H_{1:t}\}_{t=0}^\infty$ is an information state for controller~$i$. In the
  absence of this result, the information structure is not partial history
  sharing. So, we cannot identify a dynamic program for the infinite horizon
  problem. 
\end{remark}

\section {Conclusion}
\label{sec:conclusion}

Decentralized stochastic control gives rise to new conceptual challenges as
compared to centralized stochastic control. There are two solution methodologies
to overcome these challenges: (i) the person-by-person approach and (ii) the
common-information approach. The person-by-person approach provides
the structure of globally optimal control strategies and
coupled dynamic programs that determine person-by-person optimal control
strategies. The common-information approach provides the
structure of globally optimal control strategies as well as a dynamic program
that determines globally optimal control strategies. A functional optimization
problem needs to be solved to solve the dynamic program.

In practice, both the person-by-person approach and the common information
approach need to be used in tandem to solve a decentralized stochastic control
problem. For example, in the example of Sec.~\ref{sec:d-example} we first used
the person-by-person approach to simplify the information structure of the
system and then used the common-information approach to find a dynamic
programming decomposition. Neither approach could give a complete solution on
its own. A similar tandem approach has been used for simplifying specific
information structures~\cite{M:control-sharing}, real-time
communication~\cite{WalrandVaraiya:1983}, networked control
systems~\cite{MT:NCS}. 

Therefore, a general solution methodology for decentralized stochastic control
is as follows.
\begin{enumerate}
  \item Use the person-by-person approach to simplify the information structure
    of the system.
  \item Use the common-information approach on the simplified information
    structure to identify an information-state process for the system.
  \item Obtain a dynamic program corresponding to the information-state process.
  \item Either obtain an exact analytic solution of the dynamic program (as in
    the centralized case, this is possible only for very simple models), or
    obtain an approximate numerical solution of the dynamic program (as was done
    in the example above), or prove qualitative properties of optimal solution. 
\end{enumerate}
This approach is similar to the general solution approach of centralized
stochastic control, although the last step is significantly more difficult. 

%

Although we presented the common-information approach for systems with partial
history sharing information structure, the approach is applicable to all finite
horizon decentralized control problems (and extends to infinite horizon problems
under appropriate stationarity conditions). See~\cite{Nayyar:phd} for details.

\begin{acknowledgements}
  The authors are grateful to A.\ Nayyar, D.\ Teneketzis, and S.\ Y\"uksel for
  useful discussions.
\end{acknowledgements}

\end{document}